\documentclass[12pt]{amsart}
\usepackage{xcolor}
\usepackage[utf8]{inputenc}
\usepackage{graphicx}
\usepackage{amssymb,amsmath,amsgen,amsopn,amsxtra,amsbsy,amscd}

\def\({\left(}
\def\){\right)}

\def\eb{\varepsilon}

\def\rw{\rightarrow}
\def\Dx{\Delta}
\def\nx{\nabla}

\def\Om{\Omega}

\def\R {\mathbb{R}}

\def\A {{\mathcal A}}

\def\t{\tilde}

\def \p {\partial}

\def \and{\qquad\text{and}\qquad}

\def \p {\partial}
\def\Dx{\Delta}
\newcommand{\be}{\begin{equation} }
\newcommand{\ee}{\end{equation} }

\newtheorem{proposition}{Proposition}[section]
\newtheorem{theorem}[proposition]{Theorem}

\newtheorem{lemma}[proposition]{Lemma}
\theoremstyle{definition}

\newtheorem{remark}[proposition]{Remark}

\numberwithin{equation}{section}

\def \no#1#2#3 {{\bf #1} (#3), #2.}
\def \eds#1#2#3 {#1, #2, #3.}

\begin{document}
\title[Chevron pattern equations]{Global behavior of solutions to chevron pattern equations}

\author{H.~Kalantarova$^\ast$, V.~Kalantarov$^\dagger$ and O.~Vantzos$^\ddagger$}
\address{$^\ast$Department of Material Science and Engineering, Technion-Israel Institute of Technology, 32000 Haifa, Israel}
\address{$^\dagger$Department of Mathematics, Ko{\c{c}} University, Istanbul, Turkey}
\address{$^\dagger$Azerbaijan State Oil and Industry University, Baku, Azerbaijan }
\address{$^\ddagger$Lightricks Ltd., Jerusalem, Israel}


\date{\today}

\begin{abstract}
Considering a  system of equations modeling the chevron pattern
dynamics, we show that the corresponding initial boundary value problem has a unique weak solution that continuously depends on initial data, and the semigroup generated by this problem in the phase space $X^0:= L^2(\Omega)\times L^2(\Omega)$ has a global attractor. We also provide some insight to the behavior of the system, by reducing it under special assumptions to systems of ODEs, that can in turn be studied as dynamical systems.
\end{abstract}

\maketitle

\section{Introduction}\label{s0}

\begin{figure}[h]
  \centering
  \includegraphics[width=0.7\textwidth]{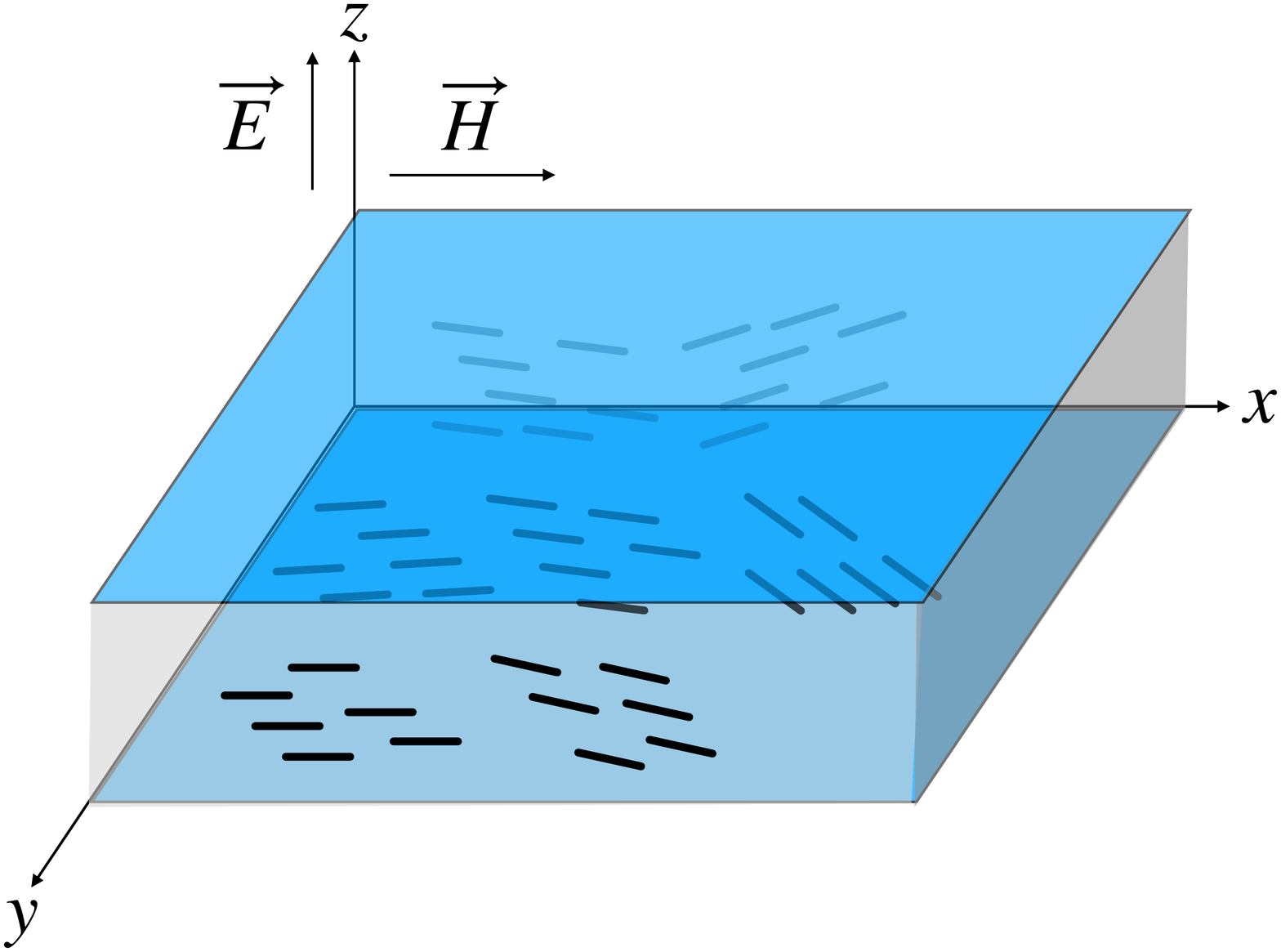} 
  \caption{\footnotesize{Cross-section of a nematic liquid crystal contained between two parallel transparent plates. The liquid is composed of rod-like particles, that are free to flow and orient themselves in 3d space, under the influence of external electric $\vec{E}$ and magnetic fields $\vec{H}$.}}
  \label{fig:cross_section}
 \end{figure}

The chevron patterns also known as the herringbone patterns in the context of the electroconvection of nematic liquid crystals, i.e.~in the electromagnetically driven motion of anisotropic liquids composed of rod-like particles, that can be oriented freely in space, were first studied by Heilmeier and Helfrich \cite{Heilmeier1970} and then in detail by Orsay group \cite{Orsay1971}. The typical experimental setup involves in general the containment of the nematic liquid between two parallel transparent plates, and the application of an AC voltage of varying frequency across the plates, often in conjunction with a magnetic field parallel to the plates, see Fig.~\ref{fig:cross_section}. Depending on the characteristics (voltage, frequency, etc.) of the external driving forces, the behavior of the fluid exhibits a wealth of nonlinear dynamical phenomena \cite{Borzsonyi1998}, \cite{Buka2012}, \cite{Demeter2000}, \cite{Kramer2001}, \cite{Scherer2000}. A common class of such phenomena feature the self-organization of the nematic liquid into convection cells, where the flow is regular and largely local to each cell. The formation of these cells is driven on a microscopical level by the interaction of the external forces with the constituent particles of the liquid, and therefore the orientation of the particles is of great importance to the dynamics. The theoretical study of this problem focuses then mainly on macroscopic models that attempt to capture directly the distribution and flow pattern (such as direction and orientation of rotation) of the convection cells, coupled with a measure of the local average orientation of the fluid particles (for instance in the form of a so-called director vector field).

\begin{figure}[h]
  \centering
  \includegraphics[width=0.7\textwidth]{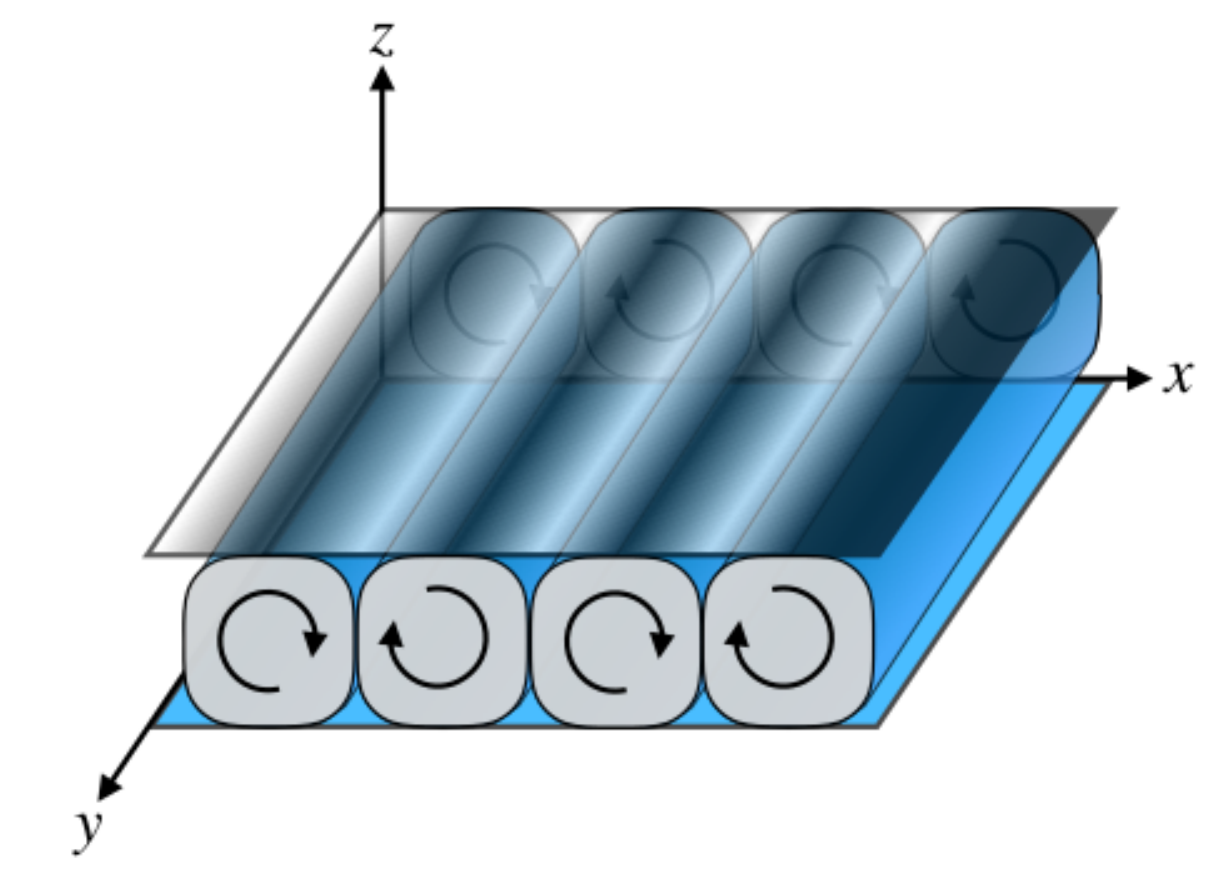} 
  \caption{\footnotesize{A typical flow pattern where the liquid self-organizes into rotating zones, called \emph{rolls}, with axes of rotation parallel to the plates and rotation orientation alternating between clockwise and counterclockwise.}}
  \label{fig:normal_roll}
 \end{figure}

\begin{figure}[h]
  \centering
  \includegraphics[width=0.7\textwidth]{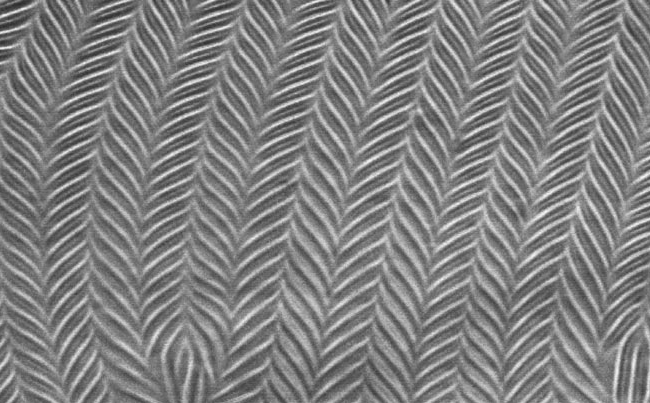} 
  \caption{\footnotesize{Experimental observation of chevron patterns, where periodic groups of rolls (observable as alternating light-dark zones) meet at an angle. Reprinted FIG. 2,(a) with permission from 
Jong-Hoon Huh, Yoshiki Hidaka, Axel G. Rossberg, and Shoichi Kai, Phys. Rev. E 61, 2769, 2000. Copyright 2000 by the
American Physical Society}}
  \label{fig:experimental}
 \end{figure}

We are interested in particular in the case where the fluid flow takes the form of \emph{rolls}, i.e. zones where the fluid rotates parallel to the plates; the rolls themselves are arranged in periodical configurations, alternating between clockwise and counter-clockwise rotation, see Fig.~\ref{fig:normal_roll}. 
The periodicity of the rolls does not hold on larger scales, which leads to interesting formations such as the titular chevron pattern, where two sequences of alternating rolls meet at an angle, see Fig.~\ref{fig:experimental}.
 The following system of equations was proposed by Rossberg et al. \cite{Rossberg_diss,Rossberg1996weakly} to model the evolution of such patterns:
\begin{eqnarray}\label{chp1}
&&\tau \partial_{t} A=A+\Dx A - \phi^2A-|A|^2A-2i c_1\phi \partial_{y} A+i\beta A\partial_{y} \phi,\\
\label{chp2}
&&\partial_{t} \phi=D_1\partial^{2}_{x}\phi+D_2\partial^{2}_{y}\phi -h\phi +\phi\lvert A\rvert^{2}-c_{2}\mbox{Im}\left[ A^{\ast}\partial_{y} A\right],
\end{eqnarray}
where $\tau$, $D_1$, $D_2$, $c_1$, $c_2$, $h$ are non-negative constants and $\beta \in \R$. The complex valued function $A$ (where $A^{\ast}$ denotes its complex conjugate) succinctly represents the phase (clockwise/ counter-clockwise), direction and amplitude (wave vector) of the periodical patterns, whereas the orientation of the liquid crystals is represented via the real valued function $\phi$, the angle of the director vector (projected in the x-y plane) with the x-axis. The parameter $\tau$ is a function of the various physical time-scales of the problem, and $D_1$ and $D_2$ are the coefficients of the anisotropic diffusion of the director field for the liquid crystal particles. The rest of the parameters reflect various interactions:  
\begin{itemize}
\item the dampening parameter $h$ measures the tendency of the director field to align with the magnetic field $\vec{H}$, corresponding to $\phi=0$,
\item the parameter $\beta$ measures the interaction between the gradient of the director field and the phase of the rolls, and
\item the parameters $c_1$ and $c_2$ control the torque that the director field and the wave vector of the rolls exert on each other; when $c_1=c_2=1$ the interaction is isotropic, but many experimentally interesting phenomena occur in the anisotropic regime.
\end{itemize}

In the literature, one can mostly find experimental \cite{Huh2000}, \cite{TaSoRo} and numerical \cite{Komineas2003}, \cite{SaMa} studies on chevron patterns, and various works on the physical derivation of the model \cite{Rossberg_diss}, \cite{Rossberg1998pattern}, \cite{Rossberg1996weakly}. It is also interesting to apply mathematical analysis tools to study this type of problem which models a physical phenomenon that features rich non-linear behavior. In our work, we show the existence, uniqueness and the continuous dependence on initial data for the weak solutions of the model. Moreover  we show the existence of a global attractor of a semigroup generated by this problem. These results are valid under the assumption that the parameter $c_{1}$ is in the range $[0,1)$.

The rest of this paper is organized as follows. In the following section, we prove well-posedness of the initial boundary value problem for the system \eqref{chp1}-\eqref{chp2} and dissipativity of the semigroup generated by this problem, provided that the coefficient of the nonlinear term $2i\phi\partial_{y} A$ in \eqref{chp1} is restricted to the range $0\leq c_{1}<1$. In  section 3, we show that the semigroup generated by the problem is a compact semigroup and has a global attractor, under the same assumption $c_{1}\in[0,1)$. In the last section, we present an argument, by reducing the model to a dynamical system, that gives insight and supports the assumed condition on the parameter $c_{1}$, by showing that $c_{1}=1$ is a critical value for the dynamics of the system.

\section*{Notation}
Throughout this paper, $\mathbb{R}$ denotes the set of real numbers, $\Omega\subset\mathbb{R}^{2}$ denotes a bounded domain with sufficiently smooth boundary denoted by $\partial \Omega$. 

\noindent $(\cdot,\cdot)$ and $\lVert\cdot\rVert$ denote the inner product and the norm induced by it in $L^{2}(\Omega)$, respectively. That is, for $f,g \in L^{2}(\Omega)$
\begin{equation}
(f,g):=\int_{\Omega}f(x,y)g^{\ast}(x,y)dxdy,\quad \lVert f\rVert:=\left(\int_{\Omega}\lvert f(x,y)\rvert^{2}dxdy\right)^{1/2}.\nonumber
\end{equation}

\section{Existence, uniqueness and dissipativity}

In this section, we study the system \eqref{chp1}-\eqref{chp2} in $\Om \subset \R^2$ under the following initial and boundary conditions
\be\label{chp3A}
A\Big|_{t=0}=A_0, \ \ \phi\Big|_{t=0}=\phi_0, \ \ A\Big|_{\p \Om}=0, \ \ \phi\Big|_{\p \Om }=0,
\ee
where $A_{0}, \phi_{0}\in L^{2}(\Omega)$ are given functions.

To prove existence and uniqueness of a weak solution of the problem we need two a priori estimates. To find these estimates, we multiply the equation \eqref{chp1} by $A^*$ and we multiply the equation \eqref{chp2} by $\phi$ and integrate them with respect to $x$ over $\Omega$. Then we take the real parts of the resulting identities and  obtain the following inequalities
\be\label{chp3}
\frac\tau 2 \frac d{dt} \| A \|^2-\|A\|^2+\| \nx A\|^2 + (\phi^2,|A|^2)+\|A\|^4_{L^4}
\le 2 c_1\lvert(\phi \p_y A,A)\rvert,
\ee
\begin{multline}\label{chp4}
\frac1 2 \frac d{dt} \| \phi\|^2+D_1\|\p_x\phi\|^2+D_2\| \p_y\phi\|^2 +(h,\phi^2)\\
\le c_2|(A^*\p_y A,\phi)|+(|\phi|^2,| A|^2),
\end{multline}
respectively, since $c_{1}, c_{2}$ are non-negative constants.
We get
\begin{eqnarray}
\label{chp3a}&&2 c_1\lvert(\phi \p_y A,A)\rvert\le c_1\|\partial_{y}A\|^2+c_1 (|A|^2, \phi^2),\\
\label{chp4a}&&c_2|(A^*\p_y A,\phi)|\le \frac{1}{2}c_2\|\p_yA\|^2+\frac{1}{2} c_2 (|A|^2, \phi^2),
\end{eqnarray}
upon application of H\"older's inequality and then Cauchy's inequality. Employing the inequality \eqref{chp3a} in \eqref{chp3} and the inequality \eqref{chp4a} in \eqref{chp4}, we obtain
\begin{multline}\label{chp5}
\frac\tau 2 \frac d{dt} \| A \|^2-\|A\|^2+(1-c_1)\| \p_yA\|^2 +\|\p_x A\|^2\\
+(1-c_1)(\phi^2,|A|^2)+\|A\|^4_{L^4(\Omega)}\le0,
\end{multline}
and
\begin{multline}\label{chp6}
\frac1 2 \frac d{dt} \| \phi\|^2+D_1\|\p_x\phi\|^2+D_2\| \p_y\phi\|^2 +(h,\phi^2)\\
\le \frac{1}{2} c_2\|\p_yA\|^2 +(1+\frac{1}{2} c_2)(|\phi|^2,| A|^2),
\end{multline}
respectively. Next, we multiply the  inequality \eqref{chp6} by a parameter $\delta>0$, whose value is to be determined later, and  add it to \eqref{chp5}
\begin{multline*}
\frac12 \frac d{dt} \left[\tau\| A \|^2 +\delta \|\phi\|^2\right]
+\delta  D_1\|\p_x\phi\|^2+\delta D_2\| \p_y\phi\|^2\\
+\delta h\|\phi\|^2-\|A\|^2+\|\p_xA\|^2+(1-c_1-\frac{1}{2}\delta c_2)\| \p_y A\|^2\\
+\big(1-c_1-\delta(1+\frac{1}{2}c_2)\big) (\phi^2,|A|^2)+\|A\|^4_{L^4(\Omega)}\le0.
\end{multline*}
Choosing
\begin{equation}
 \delta=\delta_{0}:=\frac{2(1-c_1)}{2+c_2},\quad\mbox{where}\quad c_{1}<1,\nonumber
\end{equation}
we optimize the range of value of $c_{1}$ and the number of terms with nonnegative coefficients in the above inequality, which yields
\begin{multline}\label{chp7}
\frac12 \frac d{dt} \left[\tau\| A \|^2 +\delta_{0}\|\phi\|^2\right]+ \delta_{0} D_1\|\p_x\phi\|^2+\delta_{0}D_2\| \p_y\phi\|^2\\+\delta_{0} h\|\phi\|^2
-\|A\|^2+\|\p_xA\|^2+\delta_0\| \p_y A\|^2+\|A\|^4_{L^4}\le0.
\end{multline}
Using the inequality $2\|A\|^2\le |\Om|+\|A\|^4_{L^4(\Omega)}$, where $\lvert\Omega\rvert$ denotes the measure of the domain $\Omega$, on the last term in \eqref{chp7}, we get
\begin{multline}\label{pA1}
\frac12 \frac d{dt} \left[\tau\| A \|^2 +\delta_{0}\|\phi\|^2\right] +\| A \|^2+\delta_{0} h\|\phi\|^2\\
+\delta_0\| \nabla A \|^2 +\delta_{0} D_0\|\nabla \phi\|^2\le |\Om|,
\end{multline}
where $D_0=\min\{D_1,D_2\}$. It follows from \eqref{pA1}, that
$$
 \frac d{dt} \left[\tau\| A \|^2 +\delta_{0}\|\phi\|^2\right] +k_0\left[\tau\| A \|^2 +\delta_{0}\|\phi\|^2\right] \le |\Om|,
$$
where $k_0=\min\{\tau^{-1}, h\}$.
Integrating the last inequality, we get
\be\label{chD}
\tau\| A (t)\|^2 +\delta_{0}\|\phi(t)\|^2\le [\tau\| A_0 \|^2 +\delta_{0}\|\phi_0\|^2]e^ {-k_0t}+\frac1{k_0}|\Om|.
\ee
By deriving a priori estimates \eqref{pA1} and \eqref{chD}, we have proved
\be\label{AEST1}
\|\phi(t)\|\le M_0, \ \ \|A(t)\|\le M_0, \ \ \forall t>0,
\ee
 and
\be\label{AEST2}
\int_0^T\|\nx \phi(t)\|^{2}dt \le M_T, \ \  \int_0^T\|\nx A(t)\|^{2}dt\le M_T.
\ee
Here and in what follows $M_0$ denotes a generic constant, depending only on  $\lVert A_{0}\rVert$, $\lVert\phi_{0}\rVert$, and $\lvert\Omega\rvert$; $M_{T}$ is a generic constant, which depends on $T$, $\lVert A_{0}\rVert$, $\lVert\phi_{0}\rVert$ and $\lvert\Omega\rvert$. The estimates \eqref{AEST1} and \eqref{AEST2} allow us by using the standard Galerkin method to prove the existence of a weak solution of the problem. Furthermore, employing the estimates  \eqref{AEST1} and \eqref{AEST2} we can now prove that the solution of the problem is unique. Suppose that $[A, \phi]$ and $[\t A, \t \phi]$ are solutions of the system \eqref{chp1}-\eqref{chp2} corresponding to initial data $[A_0,\phi_0]$ and $[\t A_{0}, \t \phi_{0}]$, respectively, then $[a, \Phi]:=[A-\t A, \phi-\t\phi]$ is a solution of the system
\begin{multline}\label{unc1}
\tau \p_t a =a+\Dx a-\phi^2a-(\phi^2-\t \phi^2)\t A -|A|^2A+|\t A|^2\t A\\-
2ic_1[\Phi \p_yA+\t \phi \p_y a]+i\beta[a\p_y\phi+\t A\p_y\Phi],
\end{multline}
\begin{multline}\label{unc2}
\p_t \Phi =D_1\p_x^2\Phi+D_2\p_y^2 \Phi -h\Phi+\Phi|A|^2+\tilde{\phi}(|A|^2-|\t A|^2)\\ 
-c_{2}\mbox{Im} \left[ a^{\ast}\partial_{y}A+\tilde{A}^{\ast}\partial_{y} a\right].
\end{multline}
Next we multiply the equation \eqref{unc1} by $a^*$ and the equation \eqref{unc2} by $\Phi$ in $L^2(\Omega)$ and get
\begin{multline}\label{unc3}
\frac12 \frac d{dt} \left[\tau\| a \|^2 + \|\Phi\|^2\right]-\|a\|^2+\|\nabla a\|^{2}+(\phi^{2},\lvert a\rvert^{2})\\
+ D_0\|\nx\Phi\|^2
+h\|\Phi\|^{2}\le|((\phi+\t \phi)\Phi\tilde{A}, a)|+(\lvert A\rvert\lvert\tilde{A}\rvert,\lvert a\rvert^{2})\\
+(\lvert\tilde{A}\rvert^{2},\lvert a\rvert^{2})
+2c_{1}(\lvert\Phi\rvert\lvert\partial_{y}A\rvert,\lvert a^{\ast}\rvert)
+2c_{1}(\lvert\tilde{\phi}\rvert\lvert\partial_{y}a\rvert,\lvert a^{\ast}\rvert)\\
+\beta(\lvert\tilde{A}\rvert\lvert\partial_{y}\Phi\rvert,\lvert a^{\ast}\rvert)
+(\Phi^{2},\lvert A\rvert^{2})+\lvert(\tilde{\phi}\Phi,(\lvert A\rvert+\lvert\tilde{A}\rvert)\lvert a\rvert)\rvert\\
+c_{2}\lvert(a^{\ast}\partial_{y}A,\Phi)\rvert
+c_{2}\lvert(\tilde{A}^{\ast}\partial_{y}a,\Phi)\rvert
\end{multline}
where $D_0=\min\{D_1,D_2\}$, as in \eqref{pA1}. We will estimate each term on the right hand side of the inequality \eqref{unc3}, by employing \eqref{AEST1}, H\"older's inequality, Cauchy's inequality and the Ladyzhenskaya inequality
\be\label{ineq:lad}
\|u\|_{L^4(\Om)}\le 2^{1/4}\|u\|^{1/2}\|\nx u\|^{1/2}
\ee
which is valid for each function $u\in H^1_0(\Om)$ with $\Om \subset \R^2$. We start with
\begin{eqnarray*}
 |(\phi\Phi\tilde{A}, a)|\!&\le&\!  \|\phi\|_{L^4(\Om)}\|\Phi\|_{L^4(\Om)}\|\tilde{A}\|_{L^4(\Om)}\|a\|_{L^4(\Om)}\nonumber\\
&\le&\! \frac12 \|\phi\|_{L^4(\Om)}^2\|\Phi\|_{L^4(\Om)}^2+\frac12\|\tilde{A}\|_{L^4(\Om)}^2\|a\|_{L^4(\Om)}^2\nonumber\\
&\le&\! \|\phi\|\|\nx\phi\| \|\Phi\|\|\nx\Phi\|+\|\tilde{A}\|\|\nx \tilde{A}\| \|a\|\|\nx a\|\nonumber\\
&\le&\! \eb_1 \|\nx\Phi\|^2\!+\eb_1  \|\nx a\|^2\!+\frac{M_0^2}{4\eb_1}\!\left[\|\nx\phi\|^2\|\Phi\|^2+\|\nx \tilde{A}\|^2\|a\|^2\right]
\end{eqnarray*}
which gives us the estimate of the first term on the right hand side of \eqref{unc3}
\begin{multline}\label{Ust1}
|((\phi+\t \phi)\Phi\tilde{A}, a)|\le 2\eb_1 \left[\|\nx\Phi\|^2+ \|\nx a\|^2\right]\\
+\frac{M_0^2}{4\eb_1}\left[\|\nx\phi\|^2\|\Phi\|^2+\lVert\nabla\tilde{\phi}\rVert^{2}\lVert\Phi\rVert^{2}+2\|\nx \tilde{A}\|^2\|a\|^2\right].
\end{multline}
Then we estimate the second term through the fourth terms as follows
\begin{eqnarray}\label{Ust2}
\quad\quad\quad(\lvert A\rvert\lvert\tilde{A}\rvert,\lvert a\rvert^{2})\!&\leq&\!\Vert A\rVert_{L^{4}}\lVert\tilde{A}\rVert_{L^{4}}\lVert a\rVert^{2}_{L^{4}}\\
&\leq&\!2\lVert A\rVert^{1/2}\lVert\nabla A\rVert^{1/2}\lVert \tilde{A}\rVert^{1/2}\lVert\nabla \tilde{A}\rVert^{1/2}\lVert a\rVert\lVert\nabla a\rVert\nonumber\\
&\leq&\!\lVert A\rVert \lVert\nabla A\rVert\lVert a\rVert \lVert\nabla a\rVert+\lVert \tilde{A}\rVert \lVert\nabla \tilde{A}\rVert\lVert a\rVert \lVert\nabla a\rVert\nonumber\\
&\leq&\!2\varepsilon_{2}\lVert\nabla a\rVert^{2}+\frac{M_{0}^{2}}{4\varepsilon_{2}}\left[\lVert\nabla A\rVert^{2}\lVert a\rVert^{2}+\lVert\nabla\tilde{A}\rVert^{2}\lVert a\rVert^{2}\right],\nonumber
\end{eqnarray}

\begin{eqnarray}\label{Ust3}
(\vert\tilde{A}\rvert^{2},\lvert a \rvert^{2})&\leq&\lVert\tilde{A}\rVert_{L^{4}}^{2}\lVert a\rVert_{L^{4}}^{2}\\
&\leq&2\lVert\tilde{A}\rVert\lVert\nabla\tilde{A}\rVert\lVert a\rVert\lVert\nabla a\rVert\nonumber\\
&\leq&\varepsilon_{3}\lVert\nabla a\rVert^{2}+\frac{M_{0}^{2}}{4\varepsilon_{3}}\lVert\nabla \tilde{A}\rVert^{2}\lVert a\rVert^{2},\nonumber
\end{eqnarray}

\begin{multline}\label{Ust4}
 2c_1(\lvert\Phi\rvert\lvert\partial_{y}A\rvert,\lvert a\rvert)\leq2c_{1}\lVert\Phi\rVert_{L^{4}}\lVert\nabla A\rVert\lVert a\rVert_{L^{4}}\\
\quad\quad\quad\quad\quad\quad\quad\quad\quad\quad\le 2\sqrt 2c_1 \|\Phi\|^{1/2} \|\nx \Phi\|^{1/2} \|\nabla A\|\|a\|^{1/2}\|\nx a\|^{1/2}\\
\quad\quad\quad\quad\quad\quad\ \quad\quad\quad\quad\ \le \sqrt 2 c_1 \|\Phi\|\|\nx \Phi\|\|\nabla A\|+\sqrt 2 c_1 \|\nabla A\|\|a\|\|\nx a\|\\
\quad\quad\quad\quad\quad\quad\quad\quad\quad\ \quad\ \le\! \eb_4\left[\|\nx \Phi\|^2+\|\nx a\|^2\right] +\frac{c_1^2}{2\eb_4}\|\nabla A\|^2\left[ \|a\|^2+\|\Phi\|^2\right].
\end{multline}
We can similarly estimate the other terms on the right hand side of \eqref {unc3}, and properly choosing the positive parameters $\eb_i$ in these estimates, we get the following inequality
\begin{equation}
 \frac d{dt} \left[\tau\| a \|^2 + \|\Phi\|^2\right]\\
 \le K_0\mathcal{E}(t)\left[\tau\| a \|^2 + \|\Phi\|^2\right],\nonumber
\end{equation}
where
\begin{equation}
\mathcal{E}(t)=\|\nx A(t)\|^2+ \|\nx \phi(t)\|^2+\lVert\nabla\tilde{A}(t)\rVert^{2}+\lVert\nabla\tilde{\phi}(t)\rVert^{2}.\nonumber
\end{equation}
Integrating this inequality and remembering \eqref{AEST2}, we obtain
\begin{multline}
\tau\| A(t) - \t A(t)\|^2 + \|\phi(t)-\t \phi(t)\|^2\\
\le e^ {K_0\int_0^t\mathcal{E}(\eta)d\eta}\left[\tau\| A_0-\t A_0\|^2 + \|\phi_0-\t \phi_0\|^2\right].\nonumber
\end{multline}
It follows from the last inequality that the weak solution of the problem is unique, moreover it continuously depends on initial data.\\
So we proved the following theorem

 \begin{theorem}\label{thm:diss}
The initial boundary value problem \eqref{chp1}-\eqref{chp2} and \eqref{chp3A} with $c_{1}\in[0,1)$,  has a  unique  weak solution
\begin{equation}
A,\phi \in C([0,T]; L^2(\Om)) \cap L^2([0,T]; H_0^1(\Om)),\quad \forall T>0,\nonumber  
\end{equation}
such that
\be\label{aest1}
\|\phi(t)\|\le M_0, \ \ \|A(t)\|\le M_0, \ \ \forall t>0,
\ee
 and
\be\label{aest2}
\int_0^T\|\nx \phi(t)\|^{2}dt \le M_T, \ \  \int_0^T\|\nx A(t)\|^{2}dt\le M_T, \ \forall T>0.
\ee
In other words this problem generates a continuous semigroup $S(t$), $t\ge0,$ in the phase space $X^0:= L^2(\Om)\times L^2(\Om).$
Moreover  \eqref{chD} implies that this semigroup is bounded dissipative in the phase space $X^0$.
\end{theorem}

 \begin{remark}\label{rmk:titi}(\cite{Titi2020})
 In the case $c_{1}\geq2c_{2}>0$, multiplying \eqref{chp1} and \eqref{chp2} by $c_{2}A^{\ast}$ and $2c_{1}\phi$ respectively and adding the real parts of the resulting equalities leads to the following estimate
 \begin{multline}
 \frac{d}{dt}\left(\frac{c_{1}\tau}{2}\lVert A\rVert^{2}+c_{2}\lVert\phi\rVert^{2}\right)-c_{1}\lVert A\rVert^{2}+c_{1}\lVert\nabla A\rVert^{2}+c_{1}\lVert A\rVert_{L^{4}}^{4}\\
 +2c_{2}D_{1}\lVert\partial_{x}\phi\rVert^{2}+2c_{2}D_{1}\lVert\partial_{y}\phi\rVert^{2}+2c_{2}(h,\phi^{2})\leq 0.\nonumber
 \end{multline}
 Since $\lVert A\rVert^{2}\leq\lVert A\rVert^{4}_{L^{4}}+\frac{1}{4}\lvert\Omega\rvert$ and $h>0$, we obtain the following analog of \eqref{pA1}
\begin{equation}
\frac{d}{dt}\left(\frac{c_{1}\tau}{2}\lVert A\rVert^{2}+2c_{2}\lVert\phi\rVert^{2}\right)+c_{1}\lVert\nabla A\rVert^{2}+2c_{2}d_{0}\lVert\nabla\phi\rVert^{2}
\leq\frac{1}{4}\lvert\Omega\rvert.\nonumber
\end{equation}
This inequality implies dissipativity of the system in $X^{0}$ when $c_{1}\!\geq\!2c_{2}>0$.

Thus the results of Theorem \ref{thm:diss} is valid also for $c_{1}\geq2c_{2}>0$.
 \end{remark} 

\section{Global Attractor}

In this section, we prove the existence of a global attractor for the semigroup associated with the initial boundary value problem \eqref{chp1}-\eqref{chp2} and \eqref{chp3A}. In order to establish this result, we rely on the following compactness result.

\begin{lemma}
The  semigroup $S(t): X^{0}\rightarrow X^{0}$, $t\ge0,$ is a compact semigroup.
\end{lemma}

\textbf{Proof.}
Multiplying \eqref{chp1} by $-\Delta A^{\ast}$ in $L^2$, then taking $2 Re$ parts of the obtained relation and utilizing the following inequality
\begin{multline}
-2\mbox{Re}\int_{\Omega}\lvert A\rvert^{2}A\Delta A^{\ast}dxdy
=4\int_{\Omega}\lvert\nabla A\rvert^{2}\lvert A\rvert^{2}dxdy\\
+2\mbox{Re}\int A^{2}(\nabla A^{\ast})^{2}dxdy
\geq2\int_{\Omega}\lvert\nabla A\rvert^{2}\lvert A\rvert^{2}dxdy,\nonumber
\end{multline}
we get
\begin{multline}\label{A1}
\tau\partial_{t}\lVert\nabla A\rVert^{2}-2\lVert\nabla A\rVert^{2}+2\lVert\Delta A\rVert^{2}+2(\phi^{2},\lvert\nabla A\rvert^{2})+2(\lvert\nabla A\rvert^{2},\lvert A\rvert^{2})\\
\leq 4\lvert (\phi\partial_{x}\phi,A\partial_{x} A^{\ast})\rvert+4\lvert (\phi\partial_{y}\phi,A\partial_{y} A^{\ast})\rvert\\
+4c_{1}\lvert(\phi\partial_{y}A,\Delta A^{\ast})\rvert+2\lvert\beta\rvert\lvert(A\partial_{y}\phi,\Delta A^{\ast})\rvert.
\end{multline}

We estimate each term on the right hand side of \eqref{A1} by using Cauchy's inequality, H\"older's inequality and Ladyzhenskaya inequality given in \eqref{ineq:lad}.
\begin{multline}
\label{A2a}4\lvert (A\partial_{x} A^{\ast}, \phi\partial_{x}\phi)\rvert\leq4\lVert\phi\rVert_{L^{4}}\lVert\nabla \phi\rVert_{L^{4}}\lVert A\rVert_{L^{4}}\lVert \nabla A\rVert_{L^{4}}\\
\quad\quad\quad\quad\quad\quad\ \quad\leq8\lVert\phi\rVert^{\frac{1}{2}}\lVert\nabla \phi\rVert\lVert\Delta \phi\rVert^{\frac{1}{2}}\lVert A\rVert^{\frac{1}{2}}\lVert \nabla A\rVert\lVert\Delta A\rVert^{\frac{1}{2}}\\
\quad\quad\quad\quad\quad\quad\ \quad\leq\frac{\nu_{0}}{8}\lVert\Delta\phi\rVert^{2}+\frac{32}{\nu_{0}}\lVert\phi\rVert^{2}\lVert\nabla\phi\rVert^{4}+\frac{\Vert\Delta A\rVert^{2}}{8}+32\lVert A\rVert^{2}\lVert\nabla A\rVert^{4}.
\end{multline}
Repeating exactly the same steps we estimate the second term on the right hand side of \eqref{A1} as follows
\begin{multline}\label{A2b}
4\lvert (A\partial_{y} A^{\ast}, \phi\partial_{y}\phi)\rvert\\
\leq\frac{\nu_{0}}{8}\lVert\Delta\phi\rVert^{2}+\frac{32}{\nu_{0}}\lVert\phi\rVert^{2}\lVert\nabla\phi\rVert^{4}+\frac{\Vert\Delta A\rVert^{2}}{8}+32 \lVert A\rVert^{2}\lVert\nabla A\rVert^{4}.
\end{multline}
We proceed with estimates of the third and fourth terms
\begin{eqnarray}\label{A3}
\quad4c_{1}\lvert(\phi\partial_{y}A,\Delta A)\rvert\!&\leq&\!\frac{\lVert \Delta A\rVert^{2}}{4}+16c_{1}^{2}\lVert\phi\partial_{y}A\rVert^{2}\\
&\leq&\!\frac{\lVert \Delta A\rVert^{2}}{4}+16c_{1}^{2}\lVert\phi\rVert_{L^{4}}^{2}\Vert\partial_{y}A\rVert^{2}_{L^{4}}\nonumber\\
&\leq&\!\frac{\lVert \Delta A\rVert^{2}}{4}+32c_{1}^{2}\lVert\phi\rVert\lVert\nabla\phi\rVert\lVert\nabla A\rVert\lVert\Delta A\rVert\nonumber\\
&\leq&\!\frac{\lVert \Delta A\rVert^{2}}{2}+32^{2}c_{1}^{4}\lVert\phi\rVert^{2}\lVert\nabla\phi\rVert^{2}\lVert\nabla A\rVert^{2},\nonumber
\end{eqnarray}

\begin{multline}\label{A4} 
2\lvert\beta\rvert\lvert(A\partial_{y}\phi,\Delta A)\rvert\leq\frac{\lVert\Delta A\rVert^{2}}{4}+4\beta^{2}(\lvert A\rvert^{2},\lvert\nabla\phi\rvert^{2})\\
\quad\quad\ \quad\leq\frac{\lVert\Delta A\rVert^{2}}{4}+4\beta^{2}\lVert A\rVert^{2}_{L^{4}}\lVert\nabla\phi\rVert^{2}_{L^{4}}\\
\quad\quad\quad\ \quad\quad\quad\leq\frac{\lVert\Delta A\rVert^{2}}{4}+8\beta^{2}\lVert A\rVert\lVert \nabla A\rVert\lVert\nabla\phi\rVert \lVert\Delta\phi\rVert\\
\quad\quad\quad\quad\quad\leq\!\frac{\lVert\Delta A\rVert^{2}}{4}+\frac{\nu_{0}}{4}\lVert\Delta\phi\rVert^{2}\!+\frac{64\beta^{4}}{\nu_{0}}\lVert A\rVert^{2}\lVert\nabla A\rVert^{2}\lVert \nabla \phi\rVert^{2}.
\end{multline}
Combining \eqref{A1} with \eqref{A2a}-\eqref{A4}, we obtain
\begin{multline}\label{A5}
\tau\frac{d}{dt}\lVert\nabla A\rVert^{2}\!-2\lVert\nabla A\rVert^{2}\!+2\lVert\Delta A\rVert^{2}\!+2(\phi^{2},\lvert\nabla A\rvert^{2})+2(\lvert\nabla A\rvert^{2},\lvert A\rvert^{2})\\
\quad\quad\quad\leq\frac{\nu_{0}}{2}\lVert\Delta \phi\rVert^{2}+\lVert\Delta A\rVert^{2}+\frac{64}{\nu_{0}}\lVert\phi\rVert^{2}\lVert\nabla\phi\rVert^{4}+64\lVert A\rVert^{2}\lVert\nabla A\rVert^{4}\\
+2^{10}c_{1}^{4}\lVert\phi\rVert^{2}\lVert\nabla\phi\rVert^{2}\lVert\nabla A\rVert^{2}+\frac{64 \beta^{4}}{\nu_{0}}\lVert A\rVert^{2}\lVert\nabla A\rVert^{2}\lVert\nabla\phi\rVert^{2}.
\end{multline}

On the other hand, multiplying \eqref{chp2} by
\begin{equation}
-\mathcal{L}\phi=-D_{1}\partial_{x}^{2}\phi-D_{2}\partial_{y}^{2}\phi\nonumber
\end{equation}
and using the following inequality
\begin{equation}
\nu_{0}\lVert\Delta u\rVert^{2}\leq\lVert\mathcal{L}u\rVert^{2}\leq\nu_{1}\lVert\Delta u\rVert^{2}\nonumber
\end{equation}
which is valid for each $u\in H^{2}(\Omega)\cap H^{1}_{0}(\Omega)$, we obtain
\begin{multline}\label{A6}
\frac{1}{2}\frac{d}{dt}(D_{1}\lVert\partial_{x}\phi\rVert^{2}+D_{2}\lVert\partial_{y}\phi\rVert^{2})+\nu_{0}\lVert\Delta\phi\rVert^{2}+hD_{1}\lVert\partial_{x}\phi\rVert^{2}\\
+hD_{2}\lVert\partial_{y}^{2}\phi\rVert^{2}\leq-(\phi\lvert A\rvert^{2},\mathcal{L}\phi)+c_{2}(\mbox{Im}[A^{\ast}\partial_{y}A],\mathcal{L}\phi).
\end{multline}
We estimate both terms on the right hand side of \eqref{A6} separately.
\begin{multline}
\label{A7a}\lvert(\phi\lvert A\rvert^{2},\mathcal{L}\phi)\rvert\leq D_{1}\lvert((\partial_{x}\phi)^{2},\lvert A\rvert^{2})\rvert+D_{2}\lvert((\partial_{y}\phi)^{2},\lvert A\rvert^{2})\rvert\\
\quad\quad\quad\quad\quad\quad\quad+D_{1}\lvert(\phi\partial_{x}\phi,\partial_{x}\lvert A\rvert^{2})\rvert+D_{2}\lvert(\phi\partial_{y}\phi,\partial_{y}\lvert A\rvert^{2})\rvert\\
\leq(D_{1}+D_{2})(\lvert\nabla\phi\rvert^{2},\lvert A\rvert^{2})+2D_{1}(\lvert\phi\rvert\lvert\partial_{x}\phi\rvert,\lvert A\rvert\lvert\partial_{x}A\rvert)\quad\quad\quad\quad\quad\quad\quad\quad\quad\quad\\
 \leq(D_{1}+D_{2})(\lvert\nabla\phi\rvert^{2},\lvert A\rvert^{2})+(D_{1}+D_{2})(\phi^{2},\lvert \nabla\phi\rvert^{2})
+(D_{1}+D_{2})(\lvert A^{2}\rvert,\lvert \nabla A\rvert^{2})
\end{multline}
Estimating each term on the right hand side of \eqref{A7a} requires the same sequence of arguments
\begin{multline}\label{A7}
(D_{1}+D_{2})(\lvert\nabla\phi\rvert^{2},\lvert A\rvert^{2})\leq(D_{1}+D_{2})\lVert\nabla\phi\rVert^{2}_{L^{4}}\lVert A\rVert^{2}_{L^{4}}\\
\quad\quad\quad\quad\quad\quad\quad\ \quad\quad\quad\leq2(D_{1}+D_{2})\lVert\nabla\phi\rVert\lVert\Delta\phi\rVert\lVert A\rVert\lVert \nabla A\rVert\\
\quad\quad\quad\quad\quad\quad\quad\quad\quad\quad\quad\leq\!\varepsilon\lVert\Delta\phi\rVert^{2}+\frac{(D_{1}+D_{2})^{2}}{\varepsilon}\lVert\nabla\phi\rVert^{2}\lVert A\rVert^{2}\lVert \nabla A\rVert^{2},
\end{multline}

\begin{equation}\label{A8}
(D_{1}+D_{2})(\phi^{2},\lvert \nabla\phi\rvert^{2})\leq\varepsilon\lVert\Delta\phi\rVert^{2}+\frac{(D_{1}+D_{2})^{2}}{\varepsilon}\lVert\phi\rVert^{2}\lVert\nabla\phi\rVert^{4},
\end{equation}

\begin{equation}\label{A9}
(D_{1}+D_{2})(\vert A\rvert^{2},\lvert \nabla A\rvert^{2})\leq\varepsilon\lVert\Delta A\rVert^{2}+\frac{(D_{1}+D_{2})^{2}}{\varepsilon}\lVert A\rVert^{2}\lVert\nabla A\rVert^{4}.
\end{equation}
Next, we consider the second term on the right hand side of \eqref{A6}
\begin{multline}\label{A10}
\lvert c_{2}(\mbox{Im}[A^{\ast}\partial_{y}A],\mathcal{L}\phi)\rvert\leq\frac{\varepsilon_{1}}{\nu_{1}}\lVert\mathcal{L}\phi\rVert^{2}+\frac{\nu_{1}c_{2}^{2}}{4\varepsilon_{1}}(\lvert A\rvert^{2},\lvert\partial_{y}A\rvert^{2})\\
\leq\varepsilon_{1}\lVert\Delta\phi\rVert^{2}+\varepsilon\lVert\Delta A\rVert^{2}+\left(\frac{\nu_{1}c_{2}^{2}}{4\varepsilon_{1}}\right)^{2}\frac{1}{\varepsilon}\lVert A\rVert^{2}\lVert\nabla A\rVert^{4}.
\end{multline}

\eqref{A6} together with \eqref{A7a}-\eqref{A10} implies the following estimate
\begin{multline}\label{A11}
\frac{1}{2}\frac{d}{dt}[D_{1}\lVert\partial_{x}\phi\rVert^{2}\!+\!D_{2}\lVert\partial_{y}\phi\rVert^{2}]\!+\!\nu_{0}\lVert\Delta\phi\rVert^{2}\!+\!hD_{1}\lVert\partial_{x}\phi\rVert^{2}\!\\
+\!hD_{2}\lVert\partial_{y}\phi\rVert^{2}
\leq \!(2\varepsilon\!+\varepsilon_{1})\lVert\Delta\phi\rVert^{2}\!+2\varepsilon\lVert\Delta A\rVert^{2}\!+\frac{\nu_{1}^{2}c_{2}^{4}}{2^{4}\varepsilon_{1}^{2}\varepsilon}\lVert A\rVert^{2}\lVert\nabla A\rVert^{4}\\
+\frac{(D_{1}\!+D_{2})^{2}}{\varepsilon}(\lVert\nabla\phi\rVert^{2}\lVert A\rVert^{2}\lVert\nabla A\rVert^{2}\!+\lVert\phi\rVert^{2}\lVert\nabla\phi\rVert^{4}\!+\lVert A\rVert^{2}\lVert \nabla A\rVert^{4}).
\end{multline}

We set $\varepsilon=\frac{1}{8}$ and $\varepsilon_{1}=\frac{\nu_{0}}{4}-\frac{1}{4}$ in \eqref{A11}. Then adding \eqref{A5} and \eqref{A11}, we get
\begin{multline}\label{A12}
\frac{d}{dt}[\tau\lVert\nabla A\rVert^{2}+D_{1}\lVert\partial_{x}\phi\rVert^{2}+D_{2}\lVert\partial_{y}\phi\rVert^{2}]+\nu_{0}\lVert\Delta\phi\rVert^{2}+\frac{\lVert\Delta A\rVert^{2}}{2}\\
\leq C(\lVert\nabla\phi\rVert^{2}\!+\lVert\nabla A\rVert^{2}\!+1)(\tau\lVert\nabla A\rVert^{2}\!+D_{1}\lVert\partial_{x}\phi\rVert^{2}\!+D_{2}\lVert\partial_{y}\phi\rVert^{2})
\end{multline}
where $C$ is a constant that depends on $\nu_{0}$, $c_{1}$, $c_{2}$, $\beta$, $D_{1}$, $D_{2}$ and $M_{0}$, which is a parameter defined in \eqref{aest1}.

We multiply \eqref{A12} by $t$ and rewrite the resulting inequality as follows 
\begin{multline}\label{A13}
\frac{d}{dt}(t\mathcal{E}(t))-\mathcal{E}(t)+t\left(\nu_{0}\lVert\Delta\phi\rVert^{2}+\frac{\lVert\Delta A\rVert^{2}}{2}\right)\\
\leq C(\lVert\nabla\phi\rVert^{2}+\lVert\nabla A\rVert^{2}+1)t\mathcal{E}(t)
\end{multline}
where
\begin{equation}
\mathcal{E}(t):=\tau\lVert\nabla A\rVert^{2}+D_{1}\lVert\partial_{x}\phi\rVert^{2}+D_{2}\lVert\partial_{y}\phi\rVert^{2}.\nonumber
\end{equation}
Then we integrate \eqref{A13} with respect to $t$ over $(0,t)$
\begin{multline}
t\mathcal{E}(t)+\int_{0}^{t}\eta\left(\nu_{0}\lVert\Delta\phi(\eta)\rVert^{2}+\frac{\lVert\Delta A\rVert^{2}}{2}\right)d\eta\\
\leq\int_{0}^{t}\mathcal{E}(\eta)d\eta+C\int_{0}^{t}(\lVert\nabla\phi(\eta)\rVert^{2}+\lVert\nabla A(\eta)\rVert^{2}+1)\eta\mathcal{E}(\eta)d\eta\nonumber
\end{multline}
The desired estimate follows from the above inequality by first employing the estimate \eqref{aest2} and then the Gronwall lemma
\begin{equation}\label{A14}
t\mathcal{E}(t)\leq M_{1}(t)e^{M_{2}(t)}
\end{equation}
where $M_{1}(t)$ and $M_{2}(t)$ depend only on $\lVert A_{0}\rVert$, $\lVert\phi_{0}\rVert$ and $\lvert\Omega\rvert$ as $M_{T}$ which is introduced in \eqref{AEST2}.

Inequality \eqref{A14} implies that $\forall t>0$ and $[A_{0},\phi_{0}]\in L^{2}(\Omega)\times L^{2}(\Omega)$
\begin{equation}
S(t)[A_{0},\phi_{0}]\in H^{1}_{0}(\Omega)\times H^{1}_{0}(\Omega)\nonumber
\end{equation}
i.e., the semigroup $S(t): X^{0}\rw X^{0}$, $t\ge0$ generated by the problem \eqref{chp1}-\eqref{chp3A} is a compact semigroup. \qed

Thanks to \eqref{chD} this semigroup is also bounded dissipative. Therefore the following theorem holds true (See, e.g. \cite{Temam2013}).
\begin{theorem}\label{semigroup}If $c_{1}\in[0,1)$ or $c_{1}\!\geq2c_{2}>0$, then the semigroup  
\begin{equation}
S(t):X^{0}\rightarrow X^{0},\quad t\geq0,\nonumber
\end{equation}
generated by the problem \eqref{chp1}-\eqref{chp3A}
possesses a global attractor $\A$, i.e. a compact, invariant set that attracts uniformly each bounded set of the phase space $X^0.$ 

\noindent Moreover, $\A$ is bounded in $H^{1}(\Omega)\times H^{1}(\Omega)$.
\end{theorem}

\begin{remark}
Theorems \ref{thm:diss} and \ref{semigroup} are also valid under the condition $\lvert c_{1}\rvert<1$ although the physically relevant case requires $c_{1}$ to be nonnegative as stated.
\end{remark}

\section{Basic Dynamics}

The purpose of this section is to provide some insight to the behaviour of the system of PDEs \eqref{chp1}-\eqref{chp2}, by reducing it under special assumptions to systems of ODEs, that can in turn be studied as dynamical systems. In particular, we are interested in examining a) the possibility of pattern formation, and b) the special role that Theorem \ref{semigroup} gives to the critical value $c_1=1$.

\subsection{Polar Form}

Starting again from the equations \eqref{chp1}-\eqref{chp2}, we consider the change of variables $A=\rho e^{i \psi}$ with $\rho=\rho(x,y,t)$ and $\psi=\psi(x,y,t)$. As a preliminary, we note that
\begin{equation}
 \partial_{t} A = e^{i\psi} \partial_{t}\rho + i \rho e^{i\psi}\partial_{t}\psi
\end{equation}
(and similar for the derivatives w.r.t.~$x$ and $y$), and
\begin{equation}
 \Delta A = e^{i\psi}\left(\Delta\rho -\rho\lvert\nabla\psi\rvert^2\right) + i\rho e^{i\psi}\left(\Delta\psi + \frac{2\nabla\rho\cdot\nabla\psi}{\rho}\right).
\end{equation}
Also
\begin{multline*}
 -2i c_1 \phi\partial_{y} A + i A\partial_{y}\phi = -2ic_1\phi\left(e^{i\psi} \partial_{y}\rho + i \rho e^{i\psi}\partial_{y}\psi \right) + i \beta \rho e^{i\psi}\partial_{y}\phi\\
= e^{i\psi}\left(2c_1\rho\phi\partial_{y}\psi\right) + i\rho e^{i\psi}\left(-2c_1\phi\rho^{-1}\partial_{y}\rho + \beta\partial_{y}\phi  \right)
\end{multline*}
and
\begin{equation*}
 A^*\partial_{y} A = \rho e^{-i\psi}\left(e^{i\psi} \partial_{y}\rho + i \rho e^{i\psi}\partial_{y}\psi \right) = \rho \partial_{y} \rho + i \rho^2\partial_{y} \psi.
\end{equation*}
The system then can be written, in terms of the polar variables, as
\begin{align}
 &\tau\partial_{t} \rho = \Delta\rho -\rho\lvert\nabla\psi\rvert^2 + \rho +2c_1\rho\phi\partial_{y}\psi - \phi^2\rho - \rho^3,\\
 &\tau\partial_{t} \psi = \Delta\psi + \frac{2\nabla\rho\cdot\nabla\psi}{\rho} -2c_1\phi\rho^{-1}\partial_{y}\rho + \beta\partial_{y}\phi, \\
 &\partial_{t} \phi =  \operatorname{div}(D\nabla\phi) - h\phi +\phi\rho^2 -c_2\rho^2\partial_{y} \psi,
\end{align}
where $D=\left(\begin{smallmatrix} D_1 & 0\\ 0 & D_2\end{smallmatrix}\right)$.
\subsection{Spatially uniform dynamics}

We study first the case where all the variables are constant in space, i.e.~$\rho=\rho(t)$, $\psi=\psi(t)$ and $\phi=\phi(t)$. Equivalently, this can be thought of as the result of performing the rescaling $x\rightarrow x/\epsilon$, $y\rightarrow y/\epsilon$ and dropping the higher order terms w.r.t.~$\epsilon$. In any case, we end up with the following reduced system:
\begin{align}
 \label{red1}&\tau\partial_{t} \rho= \rho(1 - \phi^2 - \rho^2),\\
 &\tau\partial_{t} \psi = 0, \\
 \label{red2}&\partial_{t} \phi =  \phi(\rho^2-h).
\end{align}
It follows immediately that the phase $\psi$ of $A$ is decoupled from the rest of the system, and can be ignored. The dynamics of the reduced system \eqref{red1}\&\eqref{red2} depend on the value of the dampening parameter $h$:
\begin{enumerate}
  \begin{figure}
  \centering
  \includegraphics[width=.49\textwidth]{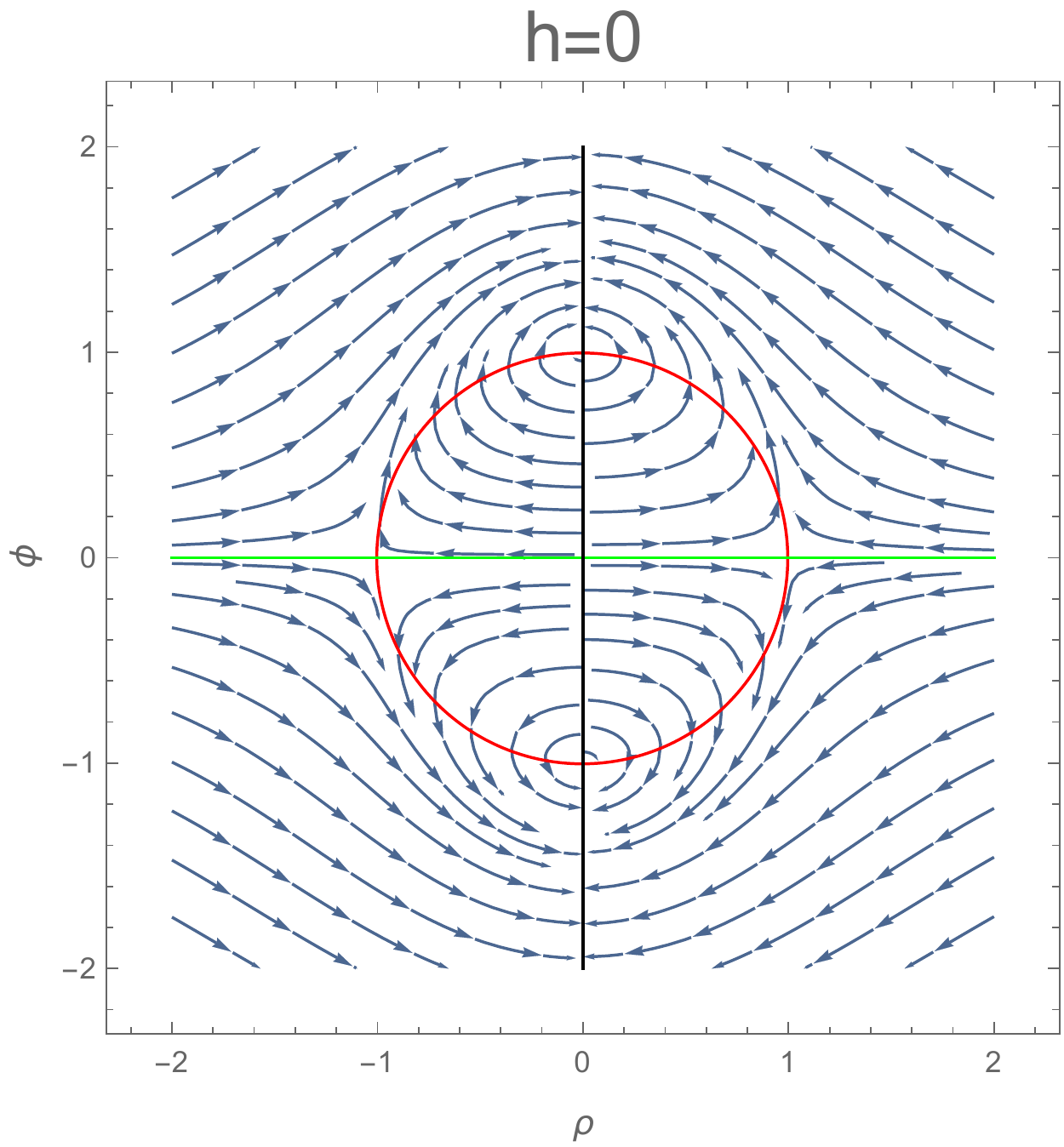} \includegraphics[width=.49\textwidth]{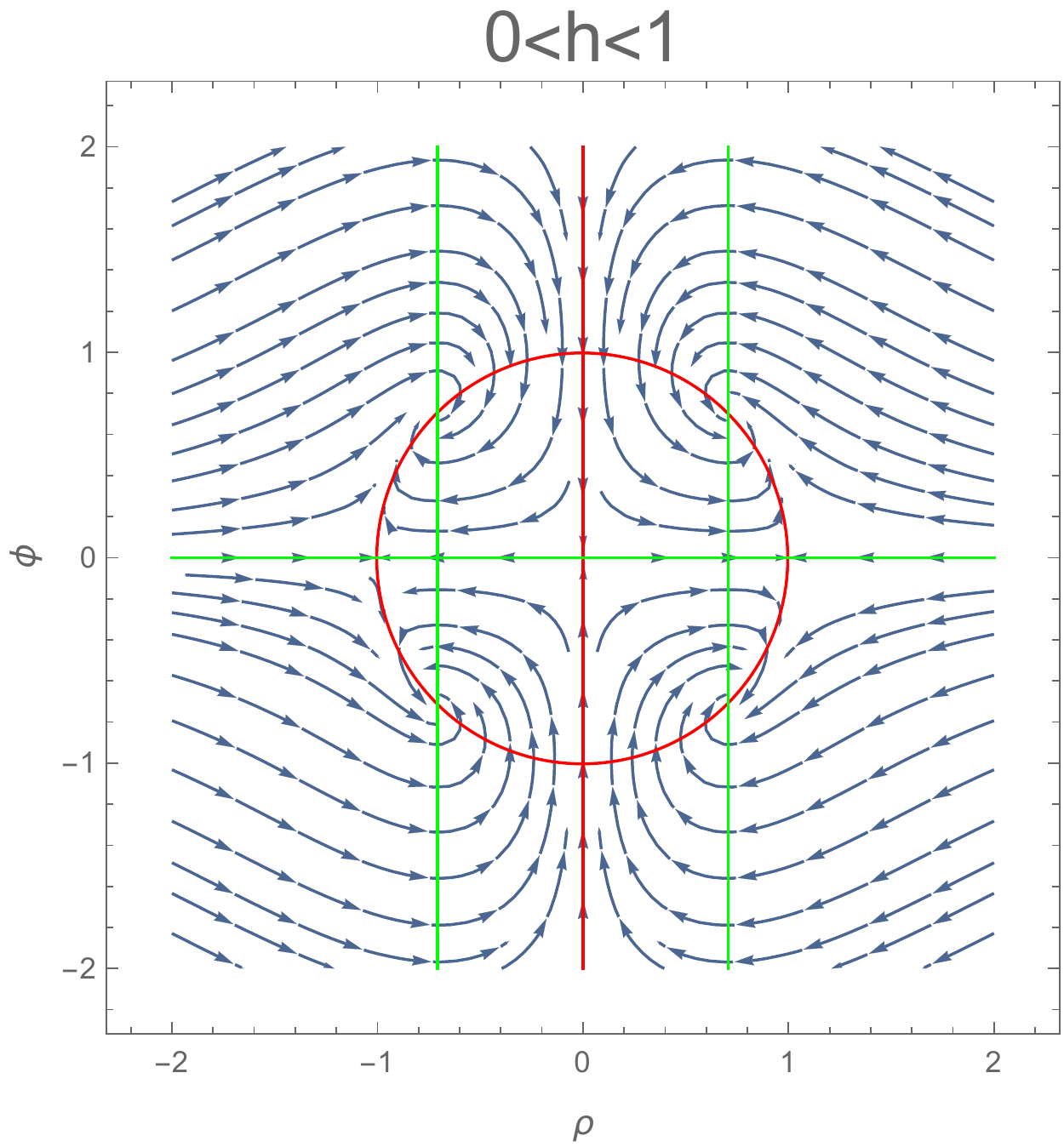}
  \caption{\footnotesize{Phase diagrams for the spatially uniform case with dampening parameter $h=0$ and dampening parameter $0<h<1$.}}
  \label{fig:h0}
 \end{figure}
 \item Fig.~\ref{fig:h0} (left), $h=0$: There are two saddle points at $(\rho,\phi)=(\pm1,0)$. The line segment $\rho=0$, $\phi\in(-1,1)$ consists of degenerate unstable critical points, whereas the rest of the $\phi$ axis is made of stable degenerate critical points.
 \item Fig.~\ref{fig:h0} (right), $0<h<1$: There are three saddle points at $(\rho,\phi)=(\pm1,0)$ and $(0,0)$, and four spiral sinks on the unit circle at $(\pm\sqrt{h},\pm\sqrt{1-h})$. Each quadrant converges to the corresponding sink.
 \item Fig.~\ref{fig:h1} (left), $h=1$: There is a saddle point at $(\rho,\phi)=(0,0)$ and two degenerate stable critical points at $(\pm 1,0)$. Points in each half-space $\phi<0$ and $\phi>0$ converge to the corresponding critical point.
 \begin{figure}
  \centering
  \includegraphics[width=.49\textwidth]{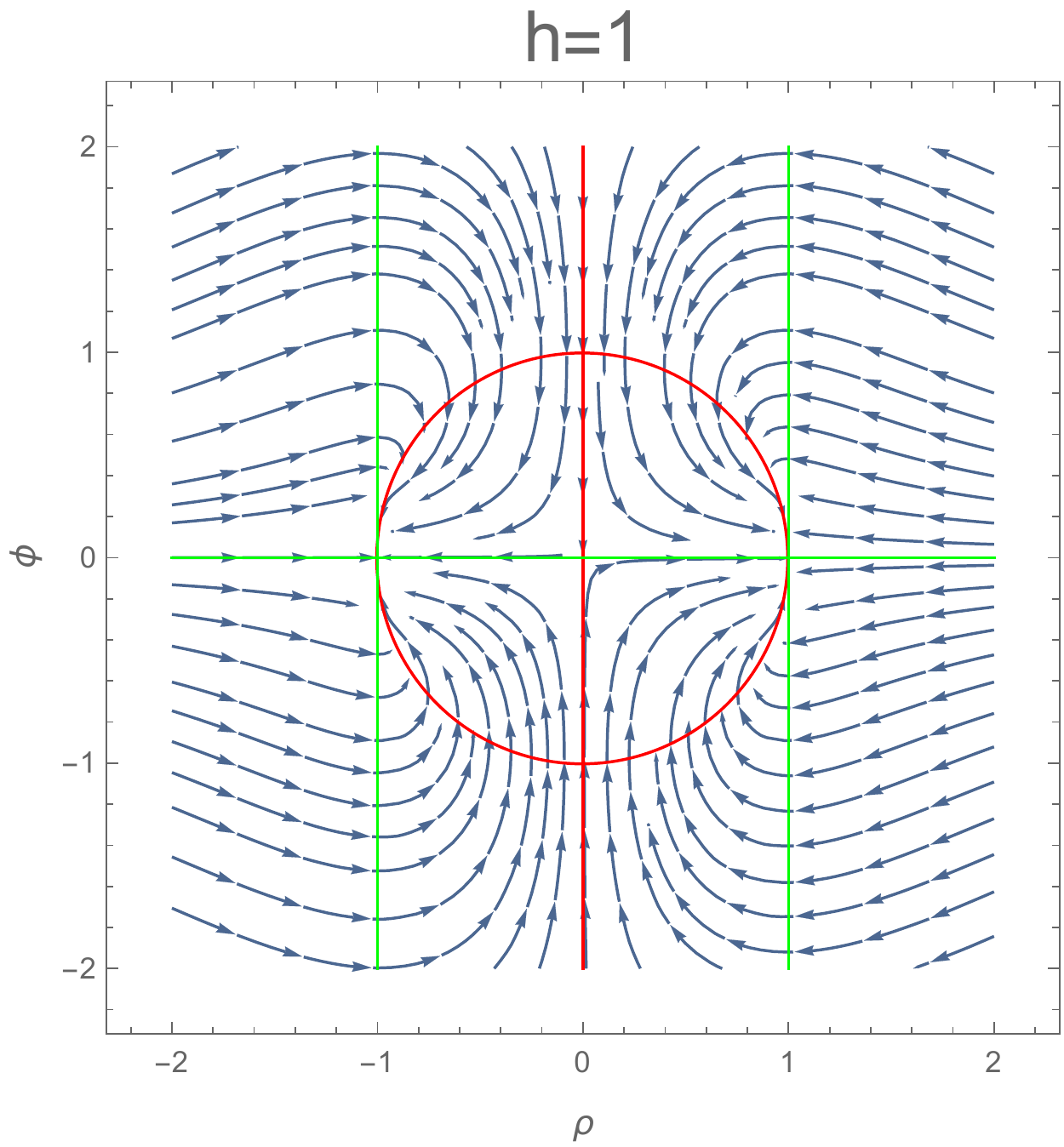}  \includegraphics[width=.49\textwidth]{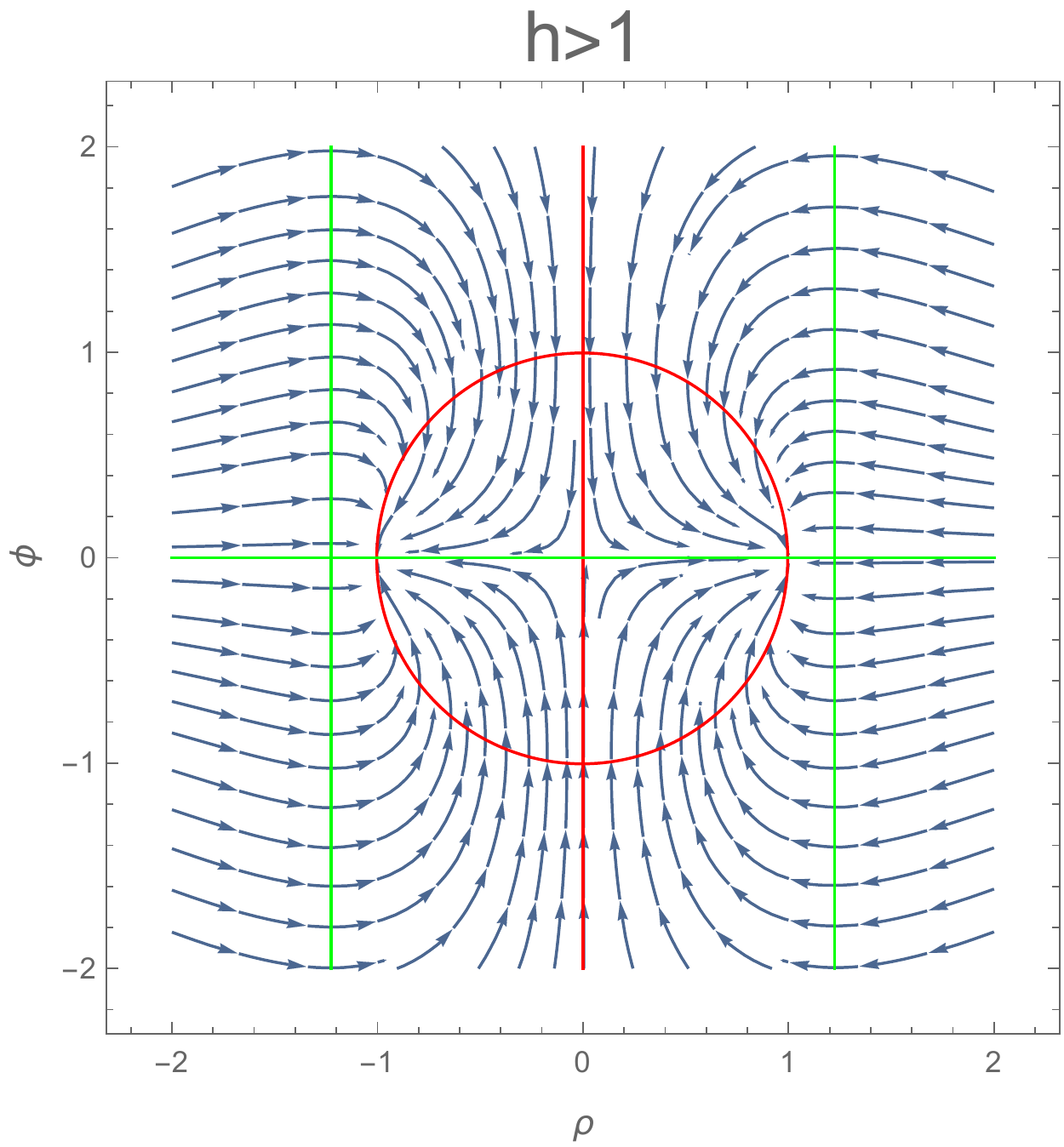}
  \caption{\footnotesize{Phase diagrams for the spatially uniform case with critical dampening parameter $h=1$ and $h>1$.}}
  \label{fig:h1}
 \end{figure}
 \item Fig.~\ref{fig:h1} (right), $h>1$: Similar behavior to the critical case $h=1$.
\end{enumerate}
Based on these cases, we expect the possibility of \emph{pattern formation} in the case when $0<h<1$, where there are two (given that $\rho$ is in fact the modulus $\lvert A\rvert$ of the complex number $A$, we are effectively only interested in the positive half-plane $\rho\geq0$) distinct non-trivial critical points  $(\rho,\phi)=(\sqrt{h},\pm\sqrt{1-h})$.

\subsection{Reduced dynamics in the presence of a phase gradient}

To reintroduce the phase variable $\psi$ into the dynamics, we assume as before that $\rho=\rho(x/\epsilon,y/\epsilon,t)$ and $\phi=\phi(x/\epsilon,y/\epsilon,t)$, but $\psi = \psi(x/\epsilon,y,t/\epsilon)$, i.e.~the variation of $\psi$ in the $y$ direction is significant. This leads to the reduced system
\begin{align*}
 &\tau\partial_{t} \rho = \rho(-\lvert\partial_{y} \psi\rvert^2 + 1 +2c_1\phi\partial_{y} \psi - \phi^2 - \rho^2),\\
 &\partial_{y}^2 \psi = 0, \\
 &\partial_{t} \phi=  \phi(\rho^2-h) - c_2\rho^2 \partial_{y} \psi.
\end{align*}
Because of the second equation the derivative $\partial_{y} \psi$ is constant. Making the change of variables $\chi = \partial_{y} \psi = const$, we end up with the following system of ODEs:
\begin{align}
 \label{eq:red_rho}&\tau\frac{d \rho}{d t} = \rho\left[(1-\rho^2) -(\phi-c_1\chi)^2 - (1-c_1^2)\chi^2\right],\\
 \label{eq:red_phi}&\frac{d \phi}{d t} =  -h\phi + \rho^2(\phi-c_2\chi).
\end{align}
The locus of the critical points of this dynamical system, for a given phase gradient $\chi$, is the set
\begin{equation}
 C_\chi = \{(\rho,\phi)\,|\,\phi=\frac{c_2\chi\rho^2}{\rho^2-h}\}
 \cap\{(\rho,\phi)\,|\,\rho^2+(\phi-c_1\chi)^2
 =1+(c_1^2-1)\chi^2\vee\rho=0\}.\nonumber
\end{equation}
The rational curve $\phi=\frac{c_2\chi\rho^2}{\rho^2-h}$ has three branches, with asymptotes at $\rho=\pm\sqrt{h}$ and $\phi=c_2\chi$, whereas the second set is the union of a circle with center $(0,c_1\chi)$ and radius $\sqrt{1+(c_1^2-1)\chi^2}$ and the line $
\rho=0$. The set always includes the point $(\rho,\phi)=(0,0)$ and, depending on the values of the parameters $c_1$ and $c_2$ and the phase gradient $\chi$, up to 3 more critical points in the positive half-plane $\rho\geq 0$, where the branches of the rational curve intersect the circle.

\begin{figure}[ht]
  \centering
  \includegraphics[width=.49\textwidth]{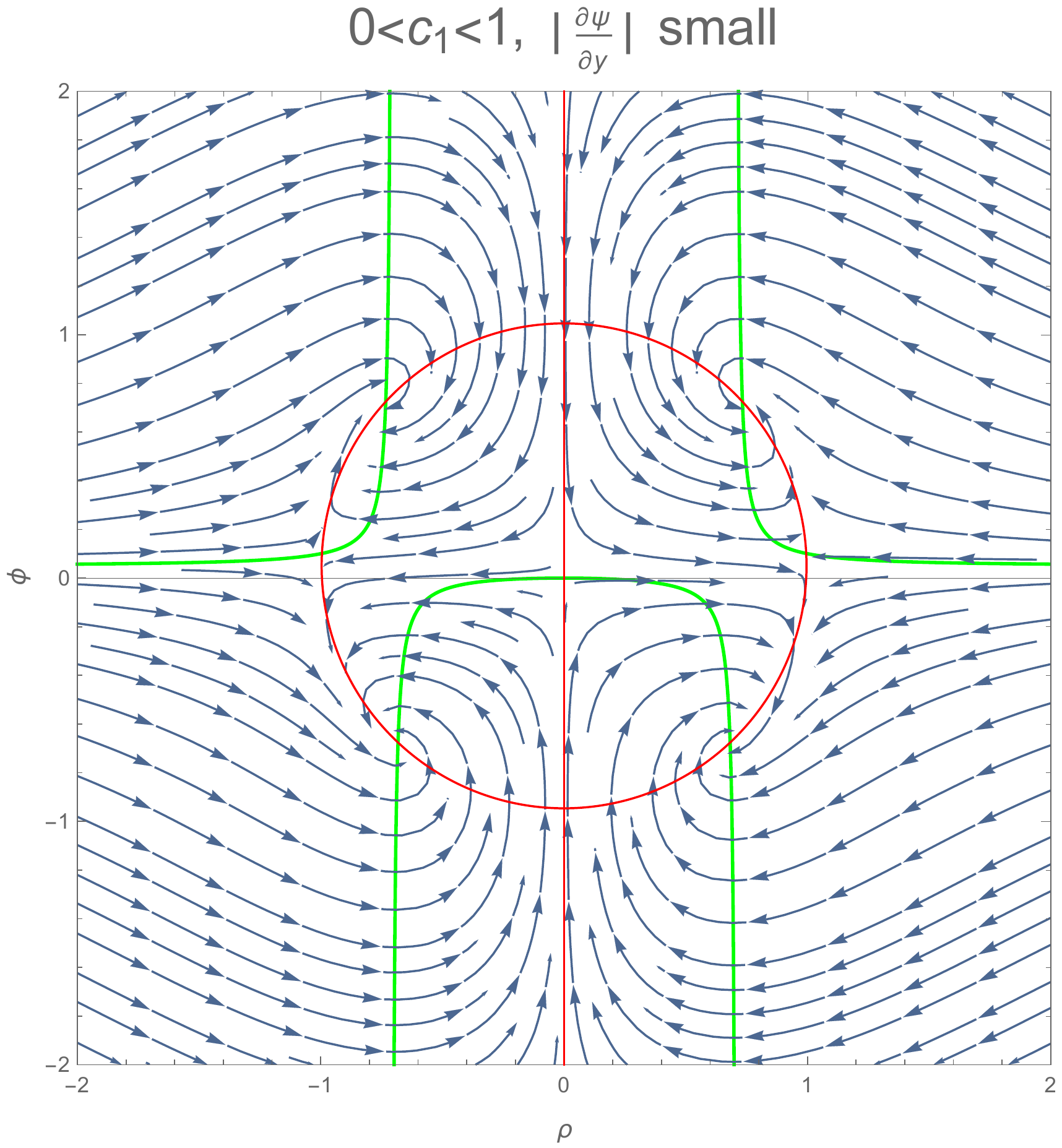} \includegraphics[width=.49\textwidth]{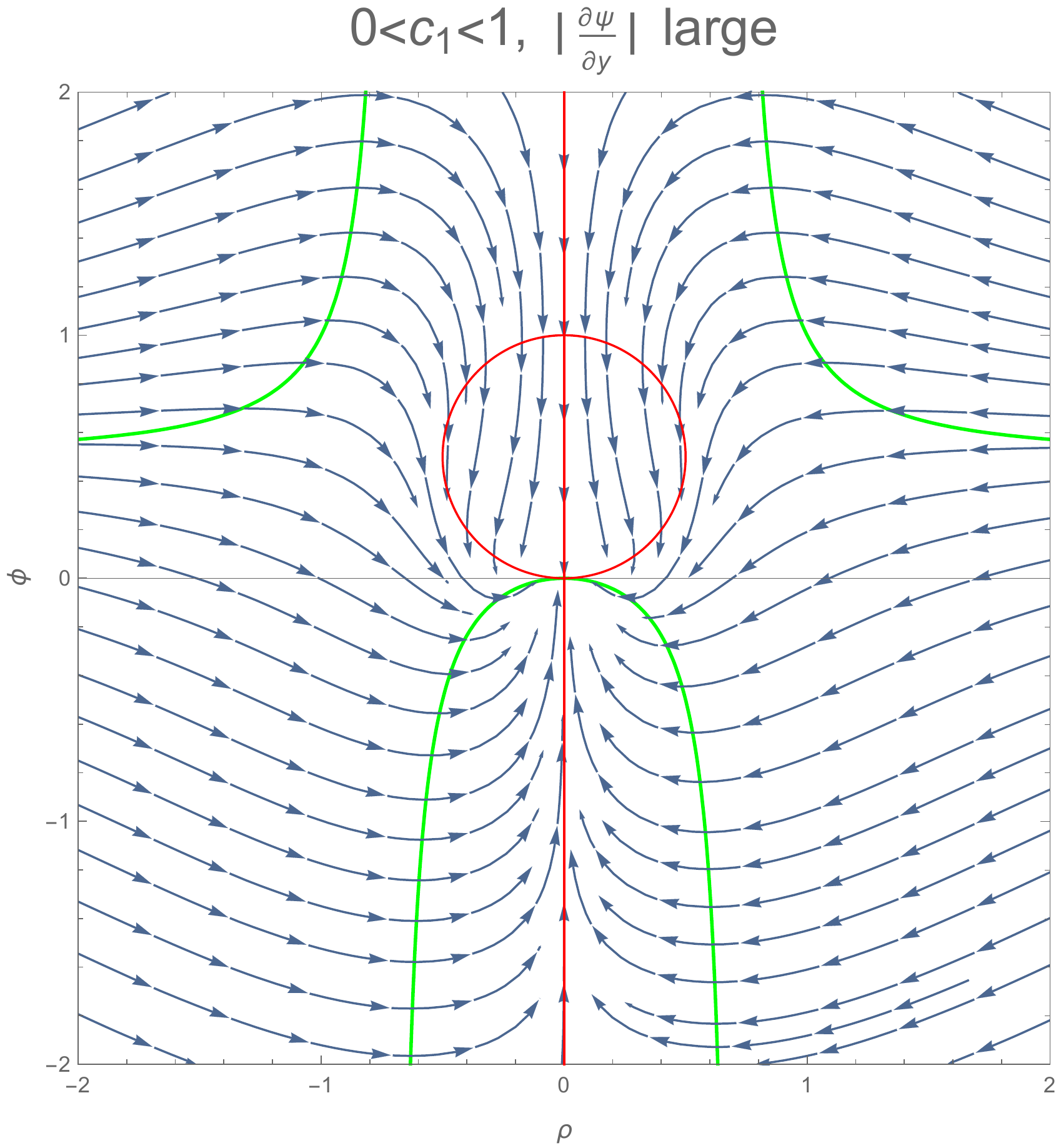}
  \caption{\footnotesize{Phase diagram in the case of $0\leq c_1<1$ (and $0<h<1$), and for different values of the phase gradient $\partial_{y}\psi$.}}
  \label{fig:c1small}
 \end{figure}

\begin{figure}[ht]
  \centering
  \includegraphics[width=.49\textwidth]{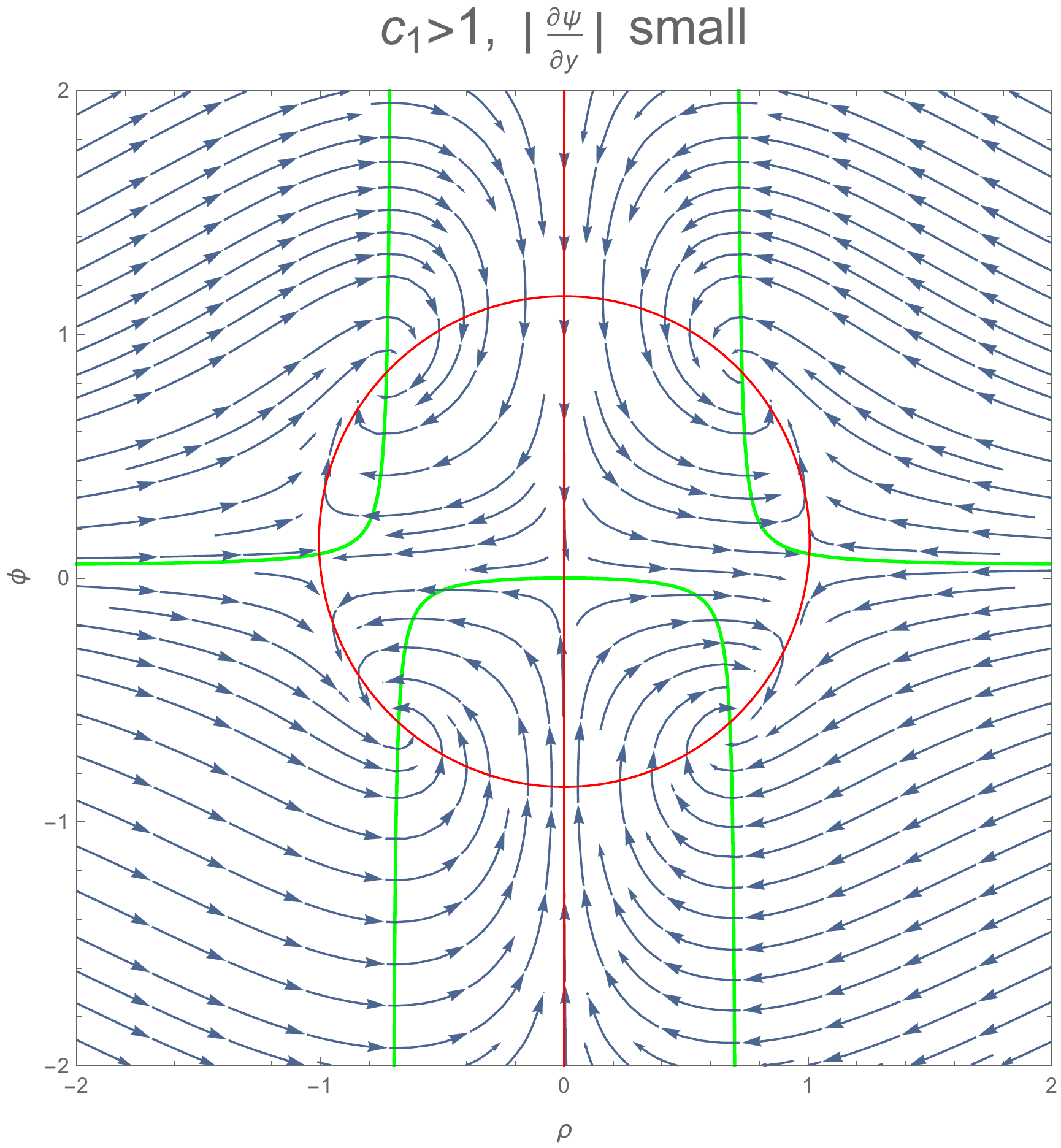} \includegraphics[width=.49\textwidth]{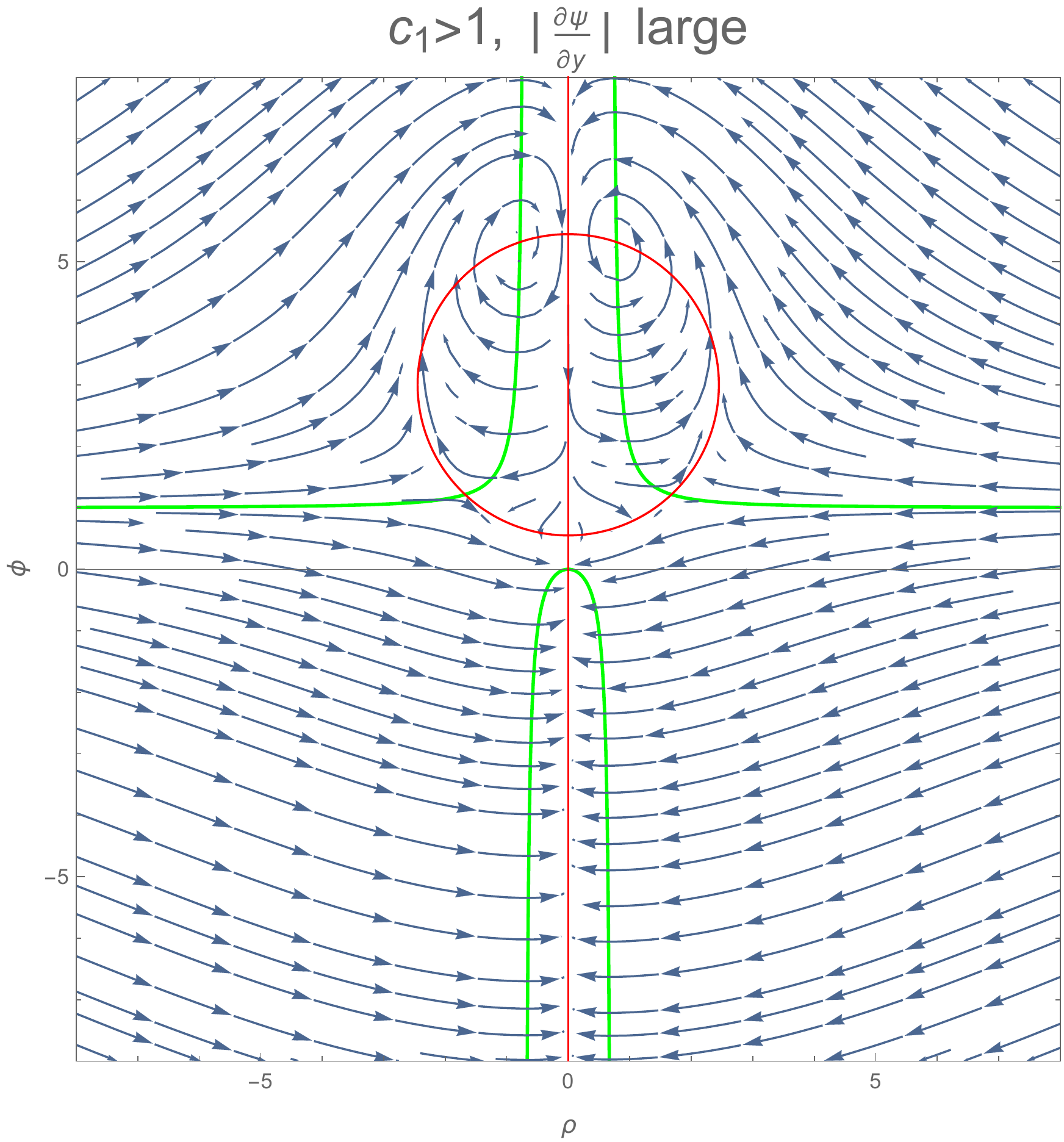}
  \caption{\footnotesize{Phase diagram in the case of $c_1>1$ (and $0<h<1$), and for different values of the phase gradient $\partial_{y}\psi$.}}
  \label{fig:c1large}
 \end{figure}

The qualitative difference in the dynamics, between the subcritical $0\leq c_1<1$ and supercritical $c_1>1$ cases, appears to be due to the location and size of the aforementioned circle on the phase diagram (see Fig.~\ref{fig:c1small} and Fig.~\ref{fig:c1large}) and the induced presence or absence of non-trivial critical points. More specifically, for small values of the phase gradient, $\lvert \partial_{t} \psi\rvert \ll 1$, the circle is centered near the origin and has radius approximately  1, yielding dynamics similar to the ones presented in the previous section. As the phase gradient $\chi$ increases in magnitude in the case $0\leq c_1<1$, the radius $\sqrt{1+(c_1^2-1)\chi^2}$ is decreasing until eventually $1+(c_1^2-1)\chi^2<0$ and there are no other critical points but the origin (Fig.~\ref{fig:c1small}), to which all orbits are attracted. On the other hand, when $c_1>1$ the radius of the circle increases with higher values of $\lvert \partial_{t} \psi\rvert$ and there is always at least one more critical point apart from the origin (Fig.~\ref{fig:c1large}).

 \begin{remark}
In Chapter 2 of \cite{Rossberg_diss}, one can find a more in depth stability analysis for solutions of the form $A=\rho e^{i(Q x+ P y)}$ for the isotropic case $c_1=c_2=1$.
\end{remark}

\begin{remark}
The results obtained in this work hold true for the case of Neumann boundary conditions as well as periodic boundary conditions, since the parameter $h$ in the equation \eqref{chp2} is positive and the equation \eqref{chp1} involves the term $\lvert A\rvert^{2} A$. In this case instead of the inequality \eqref{ineq:lad} it is necessary to use the inequality 
 \begin{equation}
 \lVert u\rVert_{L^{4}}^{2}\leq C\lVert u\rVert\lVert u\rVert_{H^{1}},\nonumber
 \end{equation}
 which is valid for each $u\in H^{1}(\Omega)$.
\end{remark}

\bibliography{chevron_patterns_B}

\providecommand{\bysame}{\leavevmode\hbox to3em{\hrulefill}\thinspace}
\providecommand{\MR}{\relax\ifhmode\unskip\space\fi MR }
\providecommand{\MRhref}[2]{%
  \href{http://www.ams.org/mathscinet-getitem?mr=#1}{#2}
}
\providecommand{\href}[2]{#2}
\begin{thebibliography}{10}

\bibitem{Borzsonyi1998}
T.~B{\"o}rzs{\"o}nyi, A.~Buka, A.~P. Krekhov, and L.~Kramer,
  \emph{\textit{Response of a homeotropic nematic liquid crystal to rectilinear
  oscillatory shear}}, Physical Review E \textbf{58} (1998), no.~6, 7419.

\bibitem{Buka2012}
A.~Buka and L.~Kramer, \emph{\textit{Pattern formation in liquid crystals}},
  Springer Science \& Business Media, 2012.

\bibitem{Demeter2000}
G.~Demeter, \emph{\textit{Complex nonlinear behavior in optically excited
  nematic liquid crystals}}, Physical Review E \textbf{61} (2000), no.~6, 6678.

\bibitem{Orsay1971}
Orsay Liquid~Crystal Group, \emph{Ac and dc regimes of the electrohydrodynamic
  instabilities in nematic liquid crystals}, Molecular Crystals and Liquid
  Crystals \textbf{12} (1971), no.~3, 251--266.

\bibitem{Heilmeier1970}
G.~H. Heilmeier and W.~Helfrich, \emph{\textit{Orientational oscillations in
  nematic liquid crystals}}, Applied Physics Letters \textbf{16} (1970), no.~4,
  155--157.

\bibitem{Huh2000}
J-H. Huh, Y.~Hidaka, A.~G. Rossberg, and S.~Kai, \emph{\textit{Pattern
  formation of chevrons in the conduction regime in homeotropically aligned
  liquid crystals}}, Physical Review E \textbf{61} (2000), no.~3, 2769.

\bibitem{Komineas2003}
S.~Komineas, H.~Zhao, and L.~Kramer, \emph{\textit{Modulated structures in
  electroconvection in nematic liquid crystals}}, Physical Review E \textbf{67}
  (2003), no.~3, 031701.

\bibitem{Kramer2001}
L.~Kramer and W.~Pesch, \emph{\textit{Electrohydrodynamic convection in
  nematics}}, EMIS DATAREVIEWS SERIES \textbf{25} (2001), 441--454.

\bibitem{Rossberg_diss}
A.~G. Rossberg, \emph{\textit{The amplitude formalism for pattern-forming
  systems with spontaneously broken isotropy and some applications}}, 1998.

\bibitem{Rossberg1996weakly}
A.~G. Rossberg, A.~Hertrich, L.~Kramer, and W.~Pesch, \emph{\textit{Weakly
  nonlinear theory of pattern-forming systems with spontaneously broken
  isotropy}}, Physical review letters \textbf{76} (1996), no.~25, 4729.

\bibitem{Rossberg1998pattern}
A.~G. Rossberg and L.~Kramer, \emph{\textit{Pattern formation from defect
  chaos—a theory of chevrons}}, Physica D: Nonlinear Phenomena \textbf{115}
  (1998), no.~1-2, 19--28.

\bibitem{SaMa}
H.~Sakaguchi and A.~Matsuda, \emph{\textit{Chevron patterns and defect lattices
  in an anisotropic model for electroconvection}}, Physica D: Nonlinear
  Phenomena \textbf{238} (2009), no.~1, 1--8.

\bibitem{Scherer2000}
M.~A. Scherer, G.~Ahlers, F.~H{\"o}rner, and I.~Rehberg,
  \emph{\textit{Deviations from linear theory for fluctuations below the
  supercritical primary bifurcation to electroconvection}}, Physical review
  letters \textbf{85} (2000), no.~18, 3754.

\bibitem{TaSoRo}
S.~Tatsumi, M.~Sano, and A.~G. Rossberg, \emph{\textit{Observation of stable
  phase jump lines in convection of a twisted nematic liquid crystal}},
  Physical Review E \textbf{73} (2006), no.~1, 011704.

\bibitem{Temam2013}
R.~Temam, \emph{\textit{Infinite-Dimensional Dynamical Systems in Mechanics and
  Physics}}, vol.~68, Springer Science \& Business Media, 2013.

\bibitem{Titi2020}
E.~S. Titi, \emph{\textit{Oral Communication}}, 2020.

\end{thebibliography}
\bibliographystyle{amsplain}

\end{document}